\documentclass{amsart}
\usepackage[latin1]{inputenc}
\usepackage{hyperref}
\makeatletter

\newfont{\gothic}{eufm10 scaled 1100}

\theoremstyle{plain}    
\newtheorem{thm}{Theorem}[section]
\numberwithin{equation}{section} 
\numberwithin{figure}{section} 
\theoremstyle{plain}    
\newtheorem{cor}[thm]{Corollary} 
\theoremstyle{plain}    
\theoremstyle{plain}    
\newtheorem{quest}[thm]{Question} 
\theoremstyle{plain}
\newtheorem{lem}[thm]{Lemma} 
\theoremstyle{plain}    
\newtheorem{prop}[thm]{Proposition} 
\theoremstyle{plain}    
\newtheorem{Def}[thm]{Definition} 
\theoremstyle{remark}
\newtheorem{rem}[thm]{Remark}
\theoremstyle{remark}

\usepackage[arrow,matrix,curve,ps]{xy}
\xyoption{dvips}
\CompileMatrices

\makeatother

\begin{document}

\title{Seshadri constants via Lelong numbers}

\date{\today}

\author{Thomas Eckl}




\address{Thomas Eckl, Mathematisches Institut, Universität
  zu K\"oln, 50931 K\"oln, Germany}

\email{thomas.eckl@math.uni-koeln.de}

\urladdr{http://www.mi.uni-koeln.de/\~{}teckl}

\maketitle

\begin{abstract}
One of Demailly's characterization of Seshadri constants on ample line bundles
works with Lelong numbers of certain positive singular hermitian metrics. 
In this note sections of multiples of the line bundle are used to produce such 
metrics and then to deduce another formula for Seshadri constants. It is 
applied to compute Seshadri constants on blown up products of curves, to 
disprove a conjectured characterization of maximal rationally connected 
quotients and to introduce a new approach to Nagata's conjecture.
\end{abstract}


\setcounter{section}{-1}

\section{Introduction}

\noindent In 1990 Demailly introduced Seshadri constants $\epsilon(L,x)$ for 
nef line bundles $L$ on projective complex manifolds $X$ \cite{Dem90}:
\[ \epsilon(L,x) := \inf_{C \ni x} \frac{L.X}{\mathrm{mult}_x C} \]
where the infimum is taken over all irreducible curves passing through $x$. 
They refine constants $\epsilon(L)$ appearing in Seshadri's ampleness criterion
\cite[Thm.1.4.13]{LazPAG1} and quantify how much of the positivity of an ample
line bundle can be localized at a given point.

\noindent These constants gained immediately a lot of interest in algebraic
geometry; for example lower bounds on Seshadri constants were used to produce
sections in adjoint bundles \cite{Laz95}. It also turned out that explicit 
calculations of Seshadri constants are difficult in almost every concrete 
situation (see for example the work of Garcia~\cite{Gar05} on ruled surfaces) 
and it is not easier to give (interesting) upper and lower bounds for them. 
From their very definition it seems easier to determine upper bounds (by 
showing that a curve with appropriate intersection number and multiplicity 
exists) than lower bounds (via the non-existence of such curves) 
\cite{EKL95, Bau99}. 

\noindent On the other hand Demailly gave two more equivalent definitions of 
Seshadri constants (\cite[Thm.6.4]{Dem90} or 
Prop.~\ref{PointSetSeshadri-prop}):
\[  
   \epsilon(L,x) = \gamma(L,x) := \sup_{\gamma \in \mathbb{R}^+} \left\{ 
   \begin{array}{l}
   \exists\ \mathrm{singular\ metric\ } h\ \mathrm{on\ } L: i\Theta_h \geq 0,\\
   x\ \mathrm{isolated\ pole\ of\ } \Theta_h,\ 
   \nu(\Theta_h,x)=\gamma
   \end{array}  \right\}, \]
where $\nu(\Theta_h,x)$ is the Lelong number of the curvature current 
$i\Theta_h$ in $x$, and
\[ \epsilon(L,x) = \sigma(L,x) := \sup_{k \in \mathbb{N}} \frac{1}{k} s(kL,x)\]
where 
\begin{eqnarray*}
s(kL,x) & := & \max_{s \in \mathbb{N}} \left\{ 
\begin{array}{c}
\mathrm{the\ sections\ in\ } H^0(X,kL)\ \mathrm{generate\ all\ } 
s-\mathrm{jets\ in\ } \\
J^s_x{kL} = \mathcal{O}_X(kL)/\hbox{\gothic m}_x^{s+1}
\end{array} \right\}.
\end{eqnarray*}

\noindent In the nef case we still have
\[ \epsilon(L,x) \geq \gamma(L,x) \geq \sigma(L,x). \]

\noindent The definitions of $\sigma(L,x)$ and $\gamma(L,x)$ allow to give 
lower bounds for Seshadri constants by constructing sections of $kL$ with 
special properties. For $\sigma(L,x)$ this is obvious whereas for $\gamma(L,x)$
we first need a better understanding of singular hermitian metrics, their
curvature currents and Lelong numbers.

\section{Lelong numbers and Seshadri constants}

\noindent
Let us repeat the relevant definitions 
(consult \cite{Dem00} for further properties and examples): A singular 
hermitian metric $h$ on a holomorphic line bundle $L$ is
given in any trivialization 
$\theta: L_{|\Omega} \stackrel{\cong}{\longrightarrow} 
 \Omega \times \mathbb{C}$ by
\[ \parallel\!\! \xi \!\!\parallel_h = |\theta(\xi)|e^{-\phi_h(x)},\ 
   x \in \Omega,\ \xi \in L_x, \]
where the function $\phi_h \in L^1_{\mathrm{loc}}(\Omega)$ is called the 
\textit{weight} of the metric $h$ (w.r.t. the trivialization $\theta$).

\noindent The curvature current $\Theta_h$ of $L$ is given by the closed 
$(1,1)$-current $\Theta_h = \frac{i}{\pi} \partial\overline{\partial} \phi_h$.
This current exists in the sense of distribution theory because of 
$\phi_h \in L^1_{\mathrm{loc}}(\Omega)$ and is independent of the chosen
trivialization.

\noindent The most important example for the rest of these notes is the 
possibly singular metric on $L$ induced by non-zero holomorphic sections
$\sigma_1, \ldots, \sigma_N$ of $L$ which is given in any trivialization 
$\theta$ by
\[ \parallel\!\! \xi \!\!\parallel_h^2 = 
   \frac{|\theta(\xi)|^2}{|\theta(\sigma_1(x))|^2 + \ldots + 
                                               |\theta(\sigma_N(x))|^2}. \]
Then the associated weight function is
\[ \phi_h(x) = \log \left( \sum_{1 \leq j \leq N} 
                           |\theta(\sigma_j(x))|^2 \right)^\frac{1}{2} \]
which is a plurisubharmonic function, so $\Theta_h$ is a (closed) positive
current. The order of logarithmic poles of a plurisubharmonic function in a
point $x \in X$ is measured by the Lelong number
\[ \nu(\phi,x) := \liminf_{z \rightarrow x} \frac{\phi(z)}{\log |z-x|}. \]
\begin{lem}
Let $\phi(z) = \log ( \sum_{1 \leq j \leq N} 
                      |\theta(\sigma_j(x))|^2)^\frac{1}{2}$ be the weight of 
the hermitian metric on $L$ induced by the holomorphic sections 
$\sigma_1, \ldots, \sigma_N$ w.r.t. a trivialization $\theta$ on an open subset
$\Omega \subset X$. Then for every point $x \in \Omega$:
\[ \nu(\phi,x) = \min_{1 \leq j \leq N} 
   \left\{ \mathrm{ord}_x \sigma_j \right\}. \]
\end{lem}
\begin{proof}
In dimension $1$ this is clear by definition. In higher dimensions restriction 
to sufficiently general lines shows at least 
\[ \nu(\phi,x) \leq \min_{1 \leq j \leq N} 
   \left\{ \mathrm{ord}_x \sigma_j \right\}. \]
To get the other inequality set $x := 0$ and add enough monomials $z^\alpha$
with $|\alpha| = \nu$  until we get
\[ \phi(z) \leq \phi^\prime(z) = \log \left( \sum_{|\alpha| = \nu} 
        |\alpha ! \cdot z^\alpha \cdot h_\alpha(z)|^2 \right)^\frac{1}{2}\]
for some holomorphic functions $h_\alpha$ defined around $0$. These 
$h_\alpha$ are bounded in a neighborhood of $0$ hence there exist constants
$C, C^\prime > 0$ such that
\[ \phi^\prime(z) \leq \log C + \log (\sum_{|\alpha| = \nu} 
        |\alpha ! \cdot z^\alpha|^2)^\frac{1}{2} \leq 
    \log C^\prime + \log (|z^\alpha|^{2\nu})^\frac{1}{2}. \]
But then
\[ \liminf_{z \rightarrow x} 
   \frac{\log (|z^\alpha|^{2\nu})^\frac{1}{2}}{\log |z|} = \nu. \]
\end{proof}

\noindent Now we are able to calculate $\gamma(L,x)$ from multiplicities of 
divisors passing through~$x$:
\begin{thm} \label{Seshadri-Lelong-thm}
Let $L$ be an ample line bundle on an $n$-dimensional projective complex 
manifold $X$. Then
\[ \epsilon(L,x) = \sup_{k;D_1, \ldots, D_n \ni x} \left\{ 
   \frac{\min_{i=1,\ldots,n}\{\mathrm{mult}_x D_i\}}{k} \right\} \]
where the supremum is taken over all divisors 
$D_1, \ldots, D_n \in |kL|$ such that $x$ is an isolated point of 
$D_1 \cap \ldots \cap D_n$. 
\end{thm}
\begin{proof}
We know that for ample line bundles
\[ \epsilon(L,x) = \sigma(L,x) = \sup_k \frac{1}{k} s(kL,x). \]
so for every $\delta > 0$ there is a $k \gg 0$ such that
\[ \epsilon(L,x) - \delta \leq \frac{s}{k} := \frac{1}{k} s(kL,x) \leq 
   \epsilon(L,x). \]
Hence $H^0(X,kL)$ generates all $s$-jets in $x$ and we can find holomorphic 
sections $f_1, \ldots, f_n$ whose $s$-jets in $x$ are the monomials $z_i^s$, 
$i=1,\ldots,n$. The weight of the associated metric on $L$ has an isolated pole
of Lelong number $\frac{s}{k}$ in $x$ and in the limit we get
\[ \epsilon(L,x) \leq \sup_{k;D_1, \ldots, D_n \ni x} \left\{ 
   \frac{\min_{i=1,\ldots,n}\{\mathrm{mult}_x D_i\}}{k} \right\}. \]
On the other hand $n$ holomorphic sections $\sigma_1, \ldots, \sigma_n$ of $kL$
(which are the $0$-divisors of $D_1, \ldots, D_n$) define a metric $h$ on $kL$.
If we multiply $\phi_h$ by $\frac{1}{k}$ we get a metric on $L$ with isolated 
pole in $x$ hence the other inequality. 
\end{proof}

\begin{rem}
An algebraic proof of the theorem may be deduced from an observation already
used in the algebraic proof of the equality
\[ \epsilon(L,x) = \sigma(L,x), \]
cf.~\cite[5.18]{LazPAG1}: Let $C \ni x$ be any irreducible and reduced curve. 
Since $|kL|$ separates $s$-jets we can find a divisor $F_k \in |kL|$ such that
$\mathrm{mult}_x(F_k) \geq s$ and $C \not\subset F_k$. Then
\begin{eqnarray*}
k(L \cdot C) & = & (F_k \cdot C) \\
 & \geq & \mathrm{mult}_x (F_k) \cdot \mathrm{mult}_x(C) \\
 & \geq & s \cdot \mathrm{mult}_x(C).
\end{eqnarray*}
\end{rem}

\section{Seshadri constants on blown up products of two curves}

\noindent We use Theorem~\ref{Seshadri-Lelong-thm} for computing Seshadri 
constants on products of two curves and some of their blow ups.
\begin{prop} \label{ProdCurve-prop}
Let $X = C_1 \times C_2$ be the product of two smooth projective curves,
$p_1: X \rightarrow C_1$, $p_2: X \rightarrow C_2$ the projections and 
$\hbox{\gothic a}$ a divisor of degree $a>0$ on $C_1$, $\hbox{\gothic b}$ a 
divisor of degree $b>0$ on $C_2$. Let 
$L = p_1^\ast \hbox{\gothic a} + p_2^\ast \hbox{\gothic b}$ be an ample divisor
on $X$ and $x$ any point on $X$. Then
\[ \epsilon(L,x) = \min (a,b). \]
\end{prop}
\begin{proof}
Let $F_i$ be the numerical equivalence class of the fibers of $p_i$, $i=1,2$.
Let $E$ be the exceptional divisor of the blow up of $X$ in $x$. Then
\[ (aF_1 + bF_2 - mE)(F_1 - E) = b-m\ ,\ 
   (aF_1 + bF_2 - mE)(F_2 - E) = a-m. \]
Since $F_1-E$, $F_2-E$ are the (numerical classes of the) strict transforms of 
the vertical resp. horizontal fiber through $x$ a line bundle of numerical 
class $aF_1 + bF_2 - mE$ is ample only if $b-m>0$, $a-m>0$ by Seshadri's 
ampleness criterion. This implies
\[ \epsilon(L,x) \leq \min (a,b). \]
On the other hand let $P_1 = p_1(x) \in C_1$ and $P_2 = p_2(x) \in C_2$. For
$a^\prime = ak$, $k \gg 0$, there exists an $A_2 \in |a^\prime P_1|$ such that
$A_2-P_1$ is not effective. Similarly for $b^\prime = bk$ there is a 
$B_2 \in |b^\prime P_2|$ such that $B_2-P_2$ is not effective. Setting 
$A_1 = a^\prime P_1$, $B_1 = b^\prime P_2$ we conclude that $x$ is an isolated 
point of $(A_1 + B_2) \cap (A_2+B_1)$ and
\[ \frac{\min(\mathrm{mult}_x(A_1+B_2),\mathrm{mult}_x(A_2+B_1))}{k} = 
   \frac{\min(ak,bk)}{k} = \min(a,b). \]
Now Seshadri constants are by definition numerical invariants. Since $A_1+B_2$
and $A_2+B_1$ belong to a line bundle with numerical class $k(aF_1 + bF_2)$ 
the theorem follows from Theorem~\ref{Seshadri-Lelong-thm}.
\end{proof}

\noindent Let $X = C_1 \times C_2$, $\hbox{\gothic a}$, $\hbox{\gothic b}$, 
$a$, $b$ and $L$ be as in the last proposition.
\begin{prop} \label{BlowupProd-prop}
Let $\pi_n: X^{(n)} \rightarrow X$ be the blow up of $X$ in $n$ points 
$x_1, \ldots, x_n$ where no two of them lie on the same horizontal or vertical 
fiber. Then for
\[ L^{(n)} := \pi_n^\ast L - \sum_{i=1}^n m_i E_i,\ 0 < m_i < \min(a,b),\ 
              \sum_{i=1}^n m_i \leq \max(a,b), \]
$E_i$ the exceptional divisor over $x_i$ we have 
\[ \epsilon(L^{(n)},y) = \min(a,b) \]
for every point $y$ not lying on the same horizontal or vertical fiber as one 
of the $x_i$.
\end{prop} 
\begin{proof}
Starting with Proposition~\ref{ProdCurve-prop} and using inductively the 
proposition together with Seshadri's ampleness criterion  we conclude that 
$L^{(n)}$ is ample. As before we compute for the blow up $X^{(n+1)}$ of 
$X^{(n)}$ in $y$ that ($E_{n+1}$ the exceptional divisor over $y$)
\[ (\pi_{n+1}^\ast L - \sum_{i=1}^n m_i E_i - mE_{n+1})(F_1-E_{n+1}) = b-m \]
and analogously for $F_2-E_{n+1}$. Consequently
\[ \epsilon(L^{(n)},y) \leq \min(a,b). \]

\noindent It is enough to show the proposition for $a = b$ since for two nef
line bundles $L,M$  it is by definition true that 
\[ \epsilon(L+M,x) \geq \epsilon(L,x). \]

\noindent Let $P_i = p_1(x_i), P = p_1(y), Q_i = p_2(x_i), Q = p_2(y)$. We 
distinguish two cases:
\begin{itemize}
\item[(1)] $a > \sum_i m_i$.
Then there is $q$ such that $|qaP-q\sum m_i P_i|$ is a base point free linear 
system on $C_1$. Consequently we have a divisor 
$P_1^\prime+ \ldots + P^\prime_{q(a-\sum m_i)}$ in this linear system such that
all $P_j^\prime \neq P$. Similarly, for $q \gg 0$ the linear system
$|qaQ-q\sum m_i Q_i|$ is base point free on $C_2$ and contains an element
$Q_1^\prime + \ldots + Q^\prime_{q(a-\sum m_i)}$ with $Q^\prime_j \neq Q$.
Setting
\begin{eqnarray*}
A_1 & = & qaP,\ \  A_2 = q \sum m_i P_i + \sum P_j^\prime, \\
B_1 & = & q \sum m_i Q_i + \sum Q_j^\prime, \\
B_2 & = & qaQ 
\end{eqnarray*}
we can apply Theorem~\ref{Seshadri-Lelong-thm} on
\[ \pi_n^\ast L^\prime - \sum m_i E_i \equiv q(\pi_n^\ast L - \sum m_i E_i) \]
where $L^\prime = \mathcal{O}(A_1+B_1) = \mathcal{O}(A_2+B_2)$ and conclude
\[ \epsilon(L,y) \geq \frac{qa}{q} = a. \]
\item[(2)] $a = \sum_i m_i$. Then for all $q \in \mathbb{N}$ there exists a 
$q^\prime$ such that
\[ |-(q^\prime qa - q^\prime)P + q^\prime q \sum m_i P_i| \ \ \mathrm{and\ \ } 
   |q^\prime q \sum m_i Q_i - (q^\prime qa - q^\prime)Q| \]
are base point free linear systems on $C_1$ and $C_2$. Consequently we have 
divisors $P_1^{\prime\prime} + \ldots + P_{q^\prime}^{\prime\prime}$  and
$Q_1^{\prime\prime} + \ldots + Q_{q^\prime}^{\prime\prime}$ in these linear 
systems such that $P_j^{\prime\prime} \neq P$ and $Q_j^{\prime\prime} \neq Q$. 
Setting
\begin{eqnarray*}
A_1 & = & (q^\prime qa - q^\prime)P + \sum P_j^{\prime\prime},\ \  
A_2 = q^\prime q \sum m_i P_i, \\
B_1 & = &  q^\prime q \sum m_i Q_i \\
B_2 & = & (q^\prime qa - q^\prime)Q + q^\prime q \sum l_i Q_i + 
          \sum Q_j^{\prime\prime}
\end{eqnarray*}
we can again apply 
Theorem~\ref{Seshadri-Lelong-thm} on
\[ \pi_n^\ast L^\prime - \sum m_i E_i \equiv 
   q^\prime q(\pi_n^\ast L - \sum m_i E_i) \]
where $L^\prime = \mathcal{O}(A_1+B_1) = \mathcal{O}(A_2+B_2)$ and conclude
\[ \epsilon(L,y) \geq \frac{q^\prime q a - q^\prime}{q^\prime q} = 
   a - \frac{1}{q}. \]
But since $q$ can be chosen arbitrarily big this implies
\[ \epsilon(L,y) \geq  a. \]
\end{itemize}  
\end{proof}

\noindent Of course we can similarly calculate Seshadri constants of (some) 
ample line bundles on the the blown up variety when some of the $x_i$ lie 
on the same horizontal or vertical fiber or on exceptional divisors. Let us 
instead use Prop.~\ref{BlowupProd-prop} to study

\section{Fibrations with rationally connected fibers}

\noindent In~\cite{Miy86} Miyaoka presented a method to construct fibrations 
with rational connected fibers (see also Shepherd-Barrons account of this work 
in~\cite{SB92} and Bogomolov-McQuillan's approach~\cite{BoMc01}, further 
explained in~\cite{KST05}):

\noindent Let $X$ be an $n$-dimensional projective complex manifold and fix an 
ample line bundle $A$ on $X$ (a polarization of $X$). Consider sufficiently 
general complete intersection curves 
$C \in |k_1 A \cap \ldots \cap k_{n-1} A|$, $k_i \gg 0$ on $X$. Miyaoka 
observed that the leaves of a foliation $\mathcal{F} \subset T_X$ (i.e.
a saturated sheaf closed under the Lie bracket) are 
rationally connected fibers of a (rational) fibration if $\mathcal{F}_{|C}$ is
ample. 

\noindent
To get such foliations we can use the properties of the 
Harder-Narasimhan filtration: For every torsion-free
sheaf $\mathcal{E}$ on $X$ there is a unique filtration, the so-called 
Harder-Narasimhan filtration of $\mathcal{E}$, 
\[ 0 = \mathcal{E}_0 \subset \mathcal{E}_1 \subset \ldots \subset 
       \mathcal{E}_k = \mathcal{E} \]
such that the $\mathcal{G}_i = \mathcal{E}_i/\mathcal{E}_{i-1}$ are semi-stable
torsion free sheaves whose slopes satisfy
\[ \mu_{\max}(\mathcal{E}) := \mu(\mathcal{G}_1) = 
   \frac{\deg(\mathcal{G}_1)}{rk(\mathcal{G}_1)} > \ldots > 
   \mu(\mathcal{G}_k) =: \mu_{\min}(\mathcal{E}). \]
Furthermore the $\mathcal{E}_i$ are saturated in $\mathcal{E}$. Finally,
by Mehta-Ramanathans theorem the restriction of the Harder-Narasimhan 
filtration of $\mathcal{E}$ on a general complete intersection curve $C$ is 
again the Harder-Narasimhan filtration of the vector bundle 
$\mathcal{E}_{|C}$. See~\cite{Sesh82,Siu87} for an introduction, proofs and 
further properties of the Harder-Narasimhan filtration and slopes. 

\noindent For our purposes we need the following proposition explicitly shown  
in~\cite{KST05}:
\begin{prop}
Let
\[ 0 = \mathcal{F}_0 \subset \mathcal{F}_1 \subset \ldots \subset 
       \mathcal{F}_k = T_X \]
be the Harder-Narasimhan filtration of the tangent bundle of $X$ with respect 
to the polarization $A$. Assume that
$\mu(F_1) > 0$ and set
\[ j := \max\{i: \mu(\mathcal{F}_i/\mathcal{F}_{i-1}) \}. \]
For any $0 < i \leq j$ the sheaf $\mathcal{F}_i$ is a foliation on $X$, and 
$\mathcal{F}_{i|C}$ is ample. If $\mathcal{F}_C \subset T_{X|C}$ is any ample 
subbundle, then $\mathcal{F}_C$ is contained in $\mathcal{F}_{j|C}$.
\end{prop}
\begin{proof}
See \cite[Prop.29,30]{KST05}.
\end{proof}

\noindent On the other hand Campana~\cite{Cam81, Cam94} and 
Koll\'ar, Miyaoka and Mori~\cite{KMM92} constructed the maximal rationally 
connected (MRC) quotient of a projective manifold $X$. If we call the sheaf 
$\mathcal{F}_j$ from the proposition above the \textit{maximal $A$-ample part} 
of $T_X$ it is tempting to ask
\begin{quest} \label{MRC-quest}
Does the MRC quotient induce a foliation $\mathcal{F}$ such that $\mathcal{F}$ 
is the maximal $A$-ample part of $T_X$?
\end{quest} 

\noindent Unfortunately the answer to this question is negative 
already on surfaces: Let $X = C \times \mathbb{P}^1$ be the ruled product 
surface over an elliptic curve $C$ with projections $p_1: X \rightarrow C$,
$p_2: X \rightarrow \mathbb{P}^1$. Let $x_1, x_2, x_3$ be 3 points on $X$ such
that no two of them lie on the same horizontal or vertical fiber. Let 
$\pi: \widehat{X} \rightarrow X$ be the blow up of $X$ in $x_1, x_2, x_3$ and 
let
\[ L = p_1^\ast \hbox{\gothic a} + p_2^\ast \hbox{\gothic b},\ 
   \deg_C \hbox{\gothic a} = 3,\ \deg_{\mathbb{P}^1} \hbox{\gothic b} = 4. \]
$L$ is an ample line bundle on $X$ and it follows from 
Prop.~\ref{BlowupProd-prop} that 
\[ L^\prime = \pi^\ast L - 2E_1 - 2E_2 -2E_3 \]
is ample, too. (Here the $E_i$ are the exceptional divisors on $\widehat{X}$
over $x_i$.) On the other hand the MRC quotient of $\widehat{X}$ is just the 
blown-up ruling of $X$ (that is a simple consequence of surface 
classification). The induced foliation $\mathcal{F}$ has tangent sheaf 
\[ T_{\mathcal{F}} = p_2^\ast \mathcal{O}(2) - E_1 - E_2 - E_3 \]
as a local computation around the blown up points show (see~\cite{Bru04}). We 
have
\[ L^\prime . T_{\mathcal{F}} = 0 \]
and hence no complete intersection curve $C \in |kL^\prime|$ intersects
$T_{\mathcal{F}}$ positively.

\noindent Of course the situation improves if we add to $L^\prime$ some fibers 
of the projection onto $C$. Hence we may change Question~\ref{MRC-quest} and
ask for the \textit{existence} of a polarization $A$ with the desired 
properties.

\section{Seshadri constants on sets of points}

\noindent Finally we turn to Seshadri constants associated to sets of points.
\begin{Def}
Let $L$ be an ample line bundle on an irreducible projective variety $X$ and 
let $x_1, \ldots, x_r$ be arbitrary points on $X$. Then we define 
\[ \epsilon(L; x_1, \ldots , x_r) = \max \{\epsilon \geq 0 : 
   \pi^\ast L - \epsilon \cdot \sum_{i=1}^r E_i\ \mathrm{is\ nef} \} \]
where $\pi: \widetilde{X} \rightarrow X$ is the blow up of $X$ in 
$x_1, \ldots, x_r$ and $E_i \subset \widetilde{X}$ is the exceptional divisor
over $x_i$.
\end{Def}

\noindent As for a single point there are other possibilities to compute
$\epsilon(L; x_1, \ldots , x_r)$. For simplicity of notation we only consider 
the surface case:
\begin{prop} \label{PointSetSeshadri-prop}
Let $L$ be an ample line bundle on a surface $X$ and $x_1, \ldots, x_r \in X$
be arbitrary points. Then $\epsilon(L; x_1, \ldots , x_r)$ equals
\[ \inf_C \frac{L.C}{\sum_{i=1}^r \mathrm{mult}_{x_i} C}, \]
where the infimum is taken over all irreducible curves $C \subset X$,
\[ \sup_k \frac{s(kL;x_1, \ldots , x_r)}{k} \]
where $s(kL;x_1, \ldots , x_r)$ is defined as the maximal $s$ such that 
$H^0(X,kL) \rightarrow 
 \bigoplus_{i=1}^r \mathcal{O}_{X,x_i} / \hbox{\gothic m}_{x_i}^{s+1}$ is 
surjective and
\[ \sup_{k;D_1,D_2} \frac{\min_j \min_{i=1,2} (\mathrm{mult}_{x_j} D_i)}{k} \]
where the supremum is taken over all pairs of divisors $D_1, D_2 \in |kL|$ such
that $x_1, \ldots, x_r$ are isolated points in $D_1 \cap D_2$.
\end{prop}
\begin{proof}
We slightly generalize the arguments in \cite[Thm.6.4]{Dem90}: The first 
equality is a direct consequence of Seshadri's ampleness criterion. Next 
suppose that there are two divisors $D_1, D_2 \in |kL|$ such
that $x_1, \ldots, x_r$ are isolated points in $D_1 \cap D_2$. If 
$f_{1j}, f_{2j}$ are local equations of $D_1, D_2$ w.r.t. a trivialization 
$\theta_j$ in a neighborhood $\Omega_j$ of $x_j$,
\[ |\xi|_h := |\theta_j(\xi)| e^{-\phi_j} \]
with weight function $\phi_j = \frac{1}{k} \log(|f_{1j}|^2 + |f_{2j}|^2)$
defines a positive singular hermitian metric on $L$ with isolated poles at the 
$x_j$ and Lelong numbers 
\[ \gamma_j := \nu(\phi_j, x_j) = 
   \min \{ \mathrm{mult}_{x_j} D_1, \mathrm{mult}_{x_j} D_2 \}. \]
Set $c(L) = \frac{i}{\pi} \partial \overline{\partial} \phi$ on $\Omega_j$ and 
let $C$ be an irreducible curve. Then
\begin{eqnarray*}
L.C & = & \int_C c(L)\ \geq\ \sum_{j=1}^r \int_{C \cap \Omega_j} 
                           \frac{i}{\pi} \partial \overline{\partial} \phi\ 
          \geq\ \sum_{j=1}^r \gamma_j \nu(C,x_j) \\
    & = & \sum_{j=1}^r \gamma_j \mathrm{mult}_{x_j} C\ \geq\  
          \frac{\min_j \min_{i=1,2} (\mathrm{mult}_{x_j} D_i)}{k} 
          \sum_{j=1}^r \mathrm{mult}_{x_j} C
\end{eqnarray*}
because the last integral is larger than the Lelong number of the integration 
current $[C]$ w.r.t. the weight $\phi_j$ and we may apply the comparison 
theorem with the ordinary Lelong number associated to the weight $\log |z-x|$
(see \cite[2.B]{Dem00} for more details). Therefore 
\[ \epsilon(L; x_1, \ldots , x_r) \geq 
   \sup_{k;D_1,D_2} \frac{\min_j \min_{i=1,2} (\mathrm{mult}_{x_j} D_i)}{k}. \]

\noindent Now fix a coordinate system $(z_1,z_2)$ on $\Omega_j$ centered at 
$x_j$. For $s := s(kL;x_1, \ldots , x_r)$ the sections in $H^0(X,kL)$ generate
all combinations of $s$-jets at all points $x_j$ and we can find two 
holomorphic sections $f_1, f_2$ whose $s$-jets are $z_1^s, z_2^s$ at every 
$x_j$. We define a global positive singular metric by
\[ |\xi|_h := |\theta_j(\xi)| e^{-\phi_j} \]
with weight function 
$\phi_j = \frac{1}{k} \log(|\theta_j(f_1)|^2 + |\theta_j(f_2)|^2)$. Then 
$\phi_j$ has isolated poles with Lelong number $\frac{s}{k}$ at $x_j$ thus
\[ \sup_{k;D_1,D_2} \frac{\min_j \min_{i=1,2} (\mathrm{mult}_{x_j} D_i)}{k} 
   \geq \frac{s}{k}. \]
Finally by using Kodaira's vanishing theorem for 
$\mathcal{O}(p\pi^\ast L - q \sum_{j=1}^r E_j)$ with
\[ \frac{q}{p} < \epsilon(L; x_1, \ldots , x_r) \]
we can show that
\[ \sup_k \frac{s(kL;x_1, \ldots , x_r)}{k} \geq 
   \epsilon(L; x_1, \ldots , x_r). \]
\end{proof}

\section{The Nagata conjecture}

\noindent A celebrated conjecture of Nagata may be seen as a statement on 
Seshadri constants associated to sets of points in $\mathbb{P}^2$: Let 
\[ \pi: \widetilde{X} = \mathrm{Bl}_r(\mathbb{P}^2) \rightarrow \mathbb{P}^2 \]
be the blow up of the projective plane in a finite set $S \subset \mathbb{P}^2$
consisting of $r$ very general points. Denote by $H$ the pull back of a line to
$\widetilde{X}$, and let $E_1, \ldots, E_r$ be the exceptional divisors over 
the $x_j$. Nagata~\cite{Nag59} conjectured in effect that the 
$\mathbb{R}$-divisor
\[ H - \sqrt{\frac{1}{r}} \sum_{j=1}^r E_j \]
is nef on $\widetilde{X}$ provided that $r \geq 9$.

\noindent It is well known that Nagata's conjecture is implied by another 
conjecture of Harbourne and Hirschowitz about spaces $\mathcal{L}_d(n^m)$ of 
plane curves of given degree $d$ and multiplicity at least $m$ at $n$ general
points \cite{Mir99, CilMir01}. This conjecture tries to detect those of the 
spaces $\mathcal{L}_d(n^m)$ which do not have the expected dimension
\[ \max (-1, \frac{d(d+3)}{2} - n \cdot \frac{m(m+1)}{2}). \]
In particular it implies that the $\mathcal{L}_d(n^m)$ do have expected 
dimension if $d \geq 3m$ \cite{CilMir00}. 

\noindent On the other hand there is a complete classification of (non) special
systems $\mathcal{L}_d(n^m)$ for $n \leq 9$:
\begin{thm}[{\cite[Thm.2.4]{CilMir00}}] \label{upto9-nonspecial-thm}
For $n \leq 9$ the special linear systems $\mathcal{L}_d(n^m)$ are
\[ \begin{array}{rcl}
   \mathcal{L}_d(2^m) & \mathrm{with} & m \leq d \leq 2m-2 \\
   \mathcal{L}_d(3^m) & \mathrm{with} & 3m/2 \leq d \leq 2m-2 \\
   \mathcal{L}_d(5^m) & \mathrm{with} & 2m \leq d \leq (5m-2)/2 \\
   \mathcal{L}_d(6^m) & \mathrm{with} & 12m/5 \leq d \leq (5m-2)/2 \\
   \mathcal{L}_d(7^m) & \mathrm{with} & 21m/8 \leq d \leq (8m-2)/3 \\
   \mathcal{L}_d(8^m) & \mathrm{with} & 48m/17 \leq d \leq (17m-2)/6.  
   \end{array} \]

\hfill $\Box$

\end{thm}

\noindent We show now that it is not necessary to prove the 
Harbourne--Hirschowitz conjecture for every triple $(d,n,m)$ to get Nagata's 
conjecture or at least a lower bound for the Seshadri constant of $H$ and $r$ 
general points on $\mathbb{P}^2$:
\begin{thm} \label{lowbd-Nagata-thm}
Let $r > 9$ be an integer and $(d_i, m_i)$ a sequence of pairs of positive 
integers such that 
$\frac{d_i^2}{m_i^2 \cdot r} 
 \stackrel{i \rightarrow \infty}{\longrightarrow} \frac{1}{a^2} \geq 1$ and the
space 
$\mathcal{L}_{d_i}(r^{m_i+1})$ has expected dimension $\geq 0$. Then 
\[ H - a \cdot \sqrt{\frac{1}{r}} \sum_{j=1}^r E_j \]
is nef on $\widetilde{X}$. In particular, Nagata's
conjecture is true for $r$ general points in $\mathbb{P}^2$, if $a = 1$.
\end{thm}

\noindent Note furthermore that the assumptions of the theorem on the spaces 
$\mathcal{L}_{d_i}(r^{m_i+1})$ follow from the Harbourne--Hirschowitz 
conjecture since $\frac{d_i}{m_i} \rightarrow \sqrt{r}$ implies $d_i \geq 3m_i$
for $i$ big enough if $r > 9$.

\vspace{0.2cm}

\noindent \textit{Proof of the Theorem.} Note first that if we define $m(d,r)$
as the maximal $m$ such that
\[ \frac{d(d+3)}{2} - r \cdot \frac{m(m+1)}{2} \geq 0 \]
we get for fixed $r$ that
\[ \lim_{d \rightarrow \infty} \frac{m(d,r)}{d} = \frac{1}{\sqrt{r}}. \]
Next we observe the following fact: Since $\mathcal{L}_{d_i}(r^1)$ has also the
expected dimension (this is the Multiplicity One Theorem in \cite{CilMir98}) 
the intermediate spaces $\mathcal{L}_{d_i}(r^{m_i})$ and
\[ \mathcal{L}_{d_i}(k^{m_i}, (r-k)^{m_i+1}),\ k = 0, \ldots , r  \]
(denoting the space of curves having multiplicity $m_i$ in $k$ points and
multiplicity $m_i+1$ in $r-k$ points) have expected dimension.

\noindent Now let 
$\pi: \widetilde{X} = \mathrm{Bl}_r(\mathbb{P}^2) \rightarrow \mathbb{P}^2$ be 
the blow up of $\mathbb{P}^2$ in $r$ general points $x_1, \ldots, x_r$. We  
choose $m_i$ points $y_{1j}, \ldots, y_{m_ij}$ on every exceptional divisor 
$E_j$ over $x_j$ and consider the space of curves $C \subset \widetilde{X}$
passing through the points $y_{kj_0}$, $k = 1, \ldots, m_i$ for some fixed 
$j_0$ such that
\begin{itemize}
\item[(i)] $\deg \pi(C) = d_i$ and
\item[(ii)] $\pi(C)$ has multiplicity $m_i+1$ in the points 
$x_1, \ldots, x_{j_0-1}, x_{j_0+1}, \ldots, x_r$.
\end{itemize}
The expected dimension of this space is
\begin{eqnarray*} 
\dim \mathcal{L}_{d_i}(1^{m_i}, (r-1)^{m_i+1}) - m_i & > & 
\dim \mathcal{L}_{d_i}(1^{m_i}, (r-1)^{m_i+1}) - m_i - 1 = \\
 & = & \mathcal{L}_{d_i}(r^{m_i+1}). 
\end{eqnarray*}
Hence there exists a curve $C \subset \mathbb{P}^2$ in this space having 
multiplicity exactly $m_i$ in $x_{j_0}$ and whose tangent cone is given by the 
points $y_{1j_0}, \ldots, y_{m_ij_0}$.

\noindent The $r$-tuples of pairwise disjoint points in $\mathbb{P}^2$ with one
marked point is parametrized by the subvariety
\[ \xymatrix{
   W \ar@{^{(}->}[r] & ((\mathbb{P}^2)^r - D) \times \mathbb{P}^2 
                       \ar[d]_{pr} \ar[r]^(.7){pr} & {\mathbb{P}^2}\\
                     & (\mathbb{P}^2)^r - D  &
   } \] 
where $D$ collects all $r$-tuples in $(\mathbb{P}^2)^r$ containing one point 
(at least) twice. A fiber of $W$ over an $r$-tuple in $(\mathbb{P}^2)^r - D$ is
just the scheme of the $r$ points in the tuple. Since the fibers of the other 
projection are Zariski-open and hence irreducible subsets of 
$(\mathbb{P}^2)^{r-1}$ the algebraic subset $W$ is irreducible. 

\noindent Consequently no Zariski-closed subset of $W$ can dominate 
$(\mathbb{P}^2)^r - D$ and we can find a tuple $(x_1, \ldots, x_r)$ such that 
for every $j_0 = 1, \ldots , r$ there are curves $C$ as described above. 
Considering these curves as divisors of degree $d_i$ and adding them gives a
divisor $D$ which has exactly multiplicity $m_i$ in every point $x_i$ and whose
tangent cone in $x_i$ is described by the $y_{ij}$, $j = 1, \ldots, m_i$.

\noindent Choosing points $y_{ij}^\prime$ such that
\[ \{y_{i1}, \ldots, y_{im_i}\} \cap \{y_{i1}^\prime, \ldots, y_{im_i}^\prime\}
   = \emptyset,\ \ i = 1, \ldots, r \]
we get in the same way another divisor $D^\prime$ such that $D \cap D^\prime$
contains $x_1, \ldots, x_r$ as isolated points and 
$\mathrm{mult}_{x_i} D = \mathrm{mult}_{x_i} D^\prime = m_i$. The two divisors 
$D$ and $D^\prime$ may be used to derive from Prop.~\ref{PointSetSeshadri-prop}
\[ \frac{m_i}{d_i} \leq \epsilon(H; x_1, \ldots, x_r). \] 
On the other hand computing the self intersection of 
$\pi^\ast H - \frac{1}{\sqrt{r}} \sum E_i$ shows that 
\[ \epsilon(H; x_1, \ldots, x_r) \leq \frac{1}{\sqrt{r}} \]
hence the last part of the theorem on the Nagata conjecture.
\hfill $\Box$

\noindent The reason behind this approach to Nagata's conjecture is that recent
attacks on the Harbourne--Hirschowitz conjecture work with recursions which 
unfortunately have 
some gaps. But it might be worth checking if around these gaps there are 
sequences of triples $(r, d_i, m_i)$ as in the theorem. We carry out this
program with the Ciliberto-Miranda recursion~\cite{CilMir98, CilMir00}. It 
relies on a deformation argument splitting linear systems on $\mathrm{P}^2$ in
another one on $\mathrm{P}^2$ and one on $\mathrm{P}^2$ blown up in one point
(which can immediately be transformed into a linear system on $\mathrm{P}^2$
again). Studying their intersection Ciliberto and Miranda were able to prove
\begin{prop}[{\cite[Cor.3.4]{CilMir98}}]
Fix $d$, $n$ and $m$. Suppose that positive integers $k$ and $b$ exist, with
$0 < k < d$ and $0 < b < n$ such that $\mathcal{L}_{d-k-1}((n-b)^m)$ and
$\mathcal{L}_d(1^{d-k+1},b^m)$ are non-empty and non-special linear systems.

\noindent Then $\mathcal{L}_d(n^m)$ is non-empty and non-special.
\end{prop}

\noindent Here, $\mathcal{L}_d(1^{d-k+1},b^m$ is the linear system of degree 
$d$ curves on $\mathbb{P}^2$ having multiplicity $d-k+1$ in one general point 
and $m$ in $b$ other general points. To deal with such quasi-homogeneous 
systems Ciliberto and Miranda used the quadratic Cremona transformation to show
\begin{prop}[{\cite[Prop.6.2]{CilMir98}}]
Let $\mathcal{L} = \mathcal{L}_d(1^{d-m},n^m)$ with $2 \leq m \leq d$. Write
$d = qm + \mu$ with $0 \leq \mu \leq m-1$ and $n = 2h + \epsilon$ with
$\epsilon \in \{0,1\}$. If $q > h$ then $\mathcal{L}$ is non-empty and 
non-special. 
\end{prop}
\begin{cor}[{\cite[Cor.6.3]{CilMir98}}]
Let $\mathcal{L} = \mathcal{L}_d(1^{d-m+1},n^m)$ and 
$\mathcal{L}^\prime = \mathcal{L}_{d-n}(1^{d-n-m+1},n^{m-1})$. Then 
$\mathcal{L}$ is non-empty and non-special if $\mathcal{L}^\prime$ is non-empty
and non-special.
\end{cor}

\noindent Combining these results with the classification of non-special 
systems for $r \leq 9$ points in Theorem~\ref{upto9-nonspecial-thm} 
S.~Barkowski~\cite{Bar06} was able to construct sequences as in 
Theorem~\ref{lowbd-Nagata-thm} recursively on the number of points $r$. Using 
the notation of Theorem~\ref{lowbd-Nagata-thm} their properties can be 
summarized as follows:
\[ \begin{array}{c||c|c|c|c} 
   r & s^2+1 & s^2+2 & s^2+3 & s^2+4 \\ \hline
   \sqrt{r}/a & s+\frac{1}{2} & s+\frac{1}{2} & s+\frac{2}{3} & s+\frac{5}{6}\\
   \end{array} \]
Unfortunately, for $s^2+k$, $5 \leq s \leq 2s+1$ the best lower bound is still
$\frac{1}{s+1}$ which can already be deduced from Nagata's solution for a 
square number $r = (s+1)^2$ of points. Even worse, it seems not be possible to
improve these lower bounds with the Ciliberto-Miranda recursion.

\bibliographystyle{alpha}

\end{document}